\documentclass[11pt]{article}

\usepackage{amsmath, amssymb, amsthm, graphicx}
\usepackage{geometry}
\title{The Geometric Refinement Transform: A Novel Continuous Transform Space}

\author{Zachary Mullaghy, M.S.\\
\textit{Independent Researcher}}

\date{\today}

\newcommand{\N}{\mathbb{N}}

\theoremstyle{definition}

\newtheorem{theorem}{Theorem}
\newtheorem{proposition}{Proposition}

\begin{document}

\maketitle
\begin{abstract}
This work seeks to introduce a novel and general class of continuous transforms derived from hierarchical Voronoi-based refinement schemes. This transform space generalizes classical methods such as wavelets and Radon transforms by introducing parameters of refinement multiplicity, dispersion, and rotation. We rigorously establish completeness, uniqueness, invertibility, closure, and stability through frame bounds of the transform on the space of functions of Bounded Variation, and define an inner product that naturally emerges from the refinement structure that exists in L2. We also characterize regions of the parameter space that recover both wavelet-like multiscale transforms and the generalized Radon transform. Finally applications in a variety of areas are introduced. Especially relevant is an exposition on entropy formulations. This transform is well behaved on geometrically challenging regions even when non convex although the behavior becomes more challenging to disambiguate. Despite at times complicated geometry the coefficient spectrum behaves well and structure can be elucidated from even challenging geometries.
\end{abstract}

\section{Introduction}
Classical transform methods, such as the Fourier and wavelet transforms, have long provided essential tools in physics, mathematics, and engineering. The Fourier transform \cite{bracewell1999fourier} is naturally suited to globally periodic or translation-invariant functions, while wavelet expansions \cite{daubechies1992ten, mallat2008wavelet} offer localized multiscale representations, especially effective for functions with compact support or singularities. These methods are deeply rooted in harmonic analysis and have spurred significant developments in signal processing, data compression, and numerical solutions of PDEs \cite{strang1996wavelets}. However, these classical methods often rely on rigid, predefined refinement structures and impose strong constraints on region geometry, convexity, and boundary conditions.

In contrast, the transform framework we introduce is inherently geometric and offers locations in the parameter space that can be self-similar. Overall this is leveraging a hierarchical Voronoi-based refinement process that allows for arbitrary region shapes, irregular boundaries, and topological flexibility. Our construction introduces three fundamental parameters: refinement multiplicity, dispersion, and rotation — enabling the transform to adapt to the natural structure of the function's domain. These parameters are left fixed with each refinement for this work but attempts to adapt the refinement technique would be a natural extension of this work. 

This work generalizes known transforms while extending their applicability. For instance, specific limits of the parameter space recover wavelet-like behavior, while rotation averaged structures approach generalized Radon transforms. Our goal is not merely to interpolate between existing techniques, but to define a new class of transforms rooted in geometric refinement and capable of encoding localized structure, multiscale hierarchy, and anisotropy simultaneously.

\paragraph{Key contributions.} In this paper, we:
\begin{itemize}
    \item Introduce a general transform framework based on hierarchical Voronoi refinement.
    \item Prove completeness, uniqueness, invertibility, and stability via measure-weighted frame bounds.
    \item Define a natural inner product derived from the geometric refinement structure.
    \item Characterize the density of piecewise continuous and BV functions in the transform space, enabling extension to all functions in \( L^2(\Omega, \mu) \).
    \item Provide geometric and physical interpretations using simplex tilings and rotational flows.
    \item Discuss applications to fluid dynamics, computer vision, quantum gravity, and biological systems.
\end{itemize}

\subsection{Mathematical Preliminaries}
Before introducing our transform, we review foundational concepts necessary for our development.

\subsubsection{Voronoi Diagrams}
Given a set of seed points $S = \{s_i \in \mathbb{R}^n\}$, the Voronoi diagram partitions space into regions $V_i$, where each region consists of all points closer to $s_i$ than any other seed:
\begin{equation}
    V_i = \{ x \in \mathbb{R}^n \mid \| x - s_i \| \leq \| x - s_j \|, \forall j \neq i \}.
\end{equation}
This structure allows for adaptive refinement and naturally aligns with self-similar function behavior \cite{aurenhammer1991voronoi}.

\subsubsection{Bounded Variation and Function Averages}
A function $f: \mathbb{R}^n \to \mathbb{R}$ is said to be of bounded variation (BV) if its total variation is finite. This ensures that our transform remains well-defined and convergent.

The average of a function $f$ over a region $R \subset \mathbb{R}^n$ is defined as:
\begin{equation}
    \text{Average}(f, R) = \frac{1}{|R|} \int_R f(x) \,dx.
\end{equation}
This forms the basis of our transform coefficients.

\subsubsection{Function Spaces and Domain of the Transform}

Throughout this paper, we work primarily with functions defined on a measurable domain \( \Omega \subseteq \mathbb{R}^n \) equipped with a positive measure \( \mu \). Our analysis involves the following nested sequence of function spaces:

\[
\mathcal{P} \subset \mathrm{BV}(\Omega) \subset L^2(\Omega, \mu),
\]

where:
\begin{itemize}
    \item \( \mathcal{P} \) denotes the space of piecewise continuous functions on \( \Omega \),
    \item \( \mathrm{BV}(\Omega) \) is the space of functions of bounded variation,
    \item \( L^2(\Omega, \mu) \) is the Hilbert space of square-integrable functions with respect to the measure \( \mu \).
\end{itemize}

We use \( \mathcal{P} \) for constructive definitions and intuition, particularly for early refinement steps and geometric averaging. Functions in \( \mathrm{BV}(\Omega) \) provide a broad and well-behaved class for which the transform is rigorously defined. Finally, via standard density arguments, we extend the transform to all of \( L^2(\Omega, \mu) \), ensuring applicability to a wide class of physical and mathematical problems.

Each level in this space hierarchy plays a distinct role in our development:
\begin{itemize}
    \item Completeness, uniqueness, and invertibility are first proven for \( \mathcal{P} \), and extended to \( \mathrm{BV} \) via convergence arguments.
    \item Frame bounds and energy estimates are formulated in \( L^2(\Omega, \mu) \), with the inner product derived from a measure-weighted coefficient structure.
    \item Stability and closure are established using norm convergence inherited from the Hilbert space structure of \( L^2 \).
\end{itemize}

\section{The Geometric Refinement Transform}

\subsection{Building the Transform}
Start with a point in $\Omega$, call this point $P_{0,0}$ where the first index indicates the refinement step, and the second index represents an enumeration for the points at that refinement step. We introduce formally the three parameters needed to detail the behavior of the transform:

\paragraph{Refinement Multiplicity. N} This parameter determines how many daughter regions are created from each parent region at every refinement level. It governs the granularity of the transform and directly influences the number of coefficients produced per scale. In Euclidean spaces, refinement multiplicity greater than the ambient dimension plus one ensures non-degenerate refinement structures.

\paragraph{Dispersion. $\delta$} Dispersion controls how far the new generator points are placed from the parent generator toward the boundary of the parent Voronoi cell. Specifically, for each child generator, we move a fraction \( \delta \in [0, 1] \) of the way along the vector from the parent to the boundary of its region. A dispersion of 0 means no movement (generators stack), while 1 places them exactly on the boundary.

\paragraph{Rotation.} Rotation introduces an angular offset to the directional placement of child generators relative to the parent configuration. It is used to break symmetry and prevent alignment artifacts from recurring across levels of refinement. This is particularly useful when the geometry being refined has high symmetry or when self-similarity would otherwise induce degeneracy in the transform structure. The rotation follows the special orthogonal group of dimension equal to the dimension of $\Omega$

For a refinement multiplicity of $N$ in 2 dimensions, place a "small" regular $N$-gon with $P_{0,0}$ at the centroid. For an angle of $0$, orient any vertex to the right. For each vertex direction $\hat{\boldsymbol{n}}_i$ such that $i \leq N$, $i \in \N$, find the smallest $\alpha$ (distance) such that:
\begin{equation}
    P_{0,0} + \cdot \hat{\boldsymbol{n}}_i \cdot \alpha = \text{(point on boundary of } \Omega \text{ or Voronoi edge)}
\end{equation}

Then create unique daughter points according to the equation:
\begin{equation}
	    P_{0,0} + \delta \cdot \alpha \cdot \hat{\boldsymbol{n}}_i = P_{1,i} \forall i \in \N,  i \leq N
\end{equation}

We then compute the Voronoi diagram of this point set. Each point generates a new convex region $\Omega_{1,i}$ such that the set union of all $\Omega_{1,i}$ is equal to $\Omega$.

\begin{equation}
\bigcup_{i=1}^{N} \Omega_{1,i} = \Omega
\end{equation}

We then refine each region $\Omega_{1,i}$ using the refinement explored above. We allow for a modulus-based indexing scheme where $\Omega_{m_1,n_1} \subseteq \Omega_{m_2,n_2}$ if and only if:

\begin{enumerate}
    \item \( m_1 > m_2 \)
    \item \( \left\lfloor \frac{n_1}{N^{(m_1 - m_2)}} \right\rfloor = n_2 \)
\end{enumerate}

This ensures that: 

\begin{equation}
\bigcup_{i=1}^{N^m} \Omega_{m,i} = \Omega 
\end{equation}

and 

\begin{equation}
\bigcup_{i=1}^{N^{m-1}} \Omega_{m,i} = \Omega_{m-1,1}
\end{equation}

Giving the first parent region union rule for order m, but also follows for k $\leq$ N we then have

\begin{equation}
\bigcup_{i=1 + ((k-1)*N^{(m-1)}}^{k*N^{m-1}} \Omega_{m,i} = \Omega_{m-1,k}
\end{equation}

This forces regions to refine each other, and also provides an indexing structure for the refinement.

\subsection{Coefficient Definition}
Given a function $f \in BV(\Omega)$ and a hierarchy of nested Voronoi regions ${ \Omega_{m,i} }$ obtained through refinement, we define the transform coefficient $c_{m,i}$ associated with region $\Omega_{m,i}$ as:
\begin{equation}
c_{m,i} = \text{Average}(f, \Omega_{m,i}) - \text{Average}(f, \Omega_{m-1,j}),
\end{equation}
where $\Omega_{m-1,j}$ is the unique parent region of $\Omega_{m,i}$ in the hierarchy.

That is,
\begin{equation}
c_{m,i} = \frac{1}{|\Omega_{m,i}|} \int_{\Omega_{m,i}} f(x) , dx - \frac{1}{|\Omega_{m-1,j}|} \int_{\Omega_{m-1,j}} f(x) , dx.
\end{equation}
This definition ensures that each coefficient represents the deviation of the local average from its parent scale, thus encoding scale-sensitive fluctuations in $f$.

\subsection{The Role of Simplex Geometry in Refinement}

A key insight into our transform space is its deep connection to the geometry of simplices. In any dimension $n$, the minimum number of refinement points required for completeness is $n+1$. These refinement points must be arranged in a way that they do not collapse into a lower-dimensional subspace, which aligns precisely with the structure of an $n$-simplex:

\begin{itemize}
    \item In \textbf{2D}, three points must not be collinear, forming a \textbf{triangle} (2-simplex).
    \item In \textbf{3D}, four points must not be coplanar, forming a \textbf{tetrahedron} (3-simplex).
    \item In \textbf{4D}, five points must not be confined to a single 3D subspace, forming a \textbf{4-simplex}.
\end{itemize}

This structure guarantees completeness by ensuring that the refinement points fully span the space at each level.

Moreover, the rotation constraints of our transform emerge naturally from the symmetry properties of these simplices. The required refinement points distribute optimally when placed according to the \textbf{generalized Thomson problem}, which minimizes electrostatic potential energy on a hypersphere. This ensures equitable distribution, non-degenerate refinement, and intrinsic rotational constraints.

Thus, our transform space is fundamentally structured as a hierarchical refinement of simplices, ensuring mathematical rigor while providing an optimal representation of functions.

\section*{Frame Bounds and the Inner Product Structure}

Let \( f \in L^2(\Omega, \mu) \), where \( \Omega \) is a measurable domain and \( \mu \) is a positive measure defined on \( \Omega \). Let \( \{R_i\} \) be a finite collection of measurable subregions used in the decomposition of \( f \), where each \( R_i \subseteq \Omega \) and \( \mu(R_i) > 0 \). The union of each region reproduces the original region \( \Omega \)

Let \( c_i \) denote a transform coefficient associated with region \( R_i \), derived through a local projection (e.g., averaging or differencing). The coefficient-based energy of the function is then given by the weighted sum:

\[
E_f := \sum_i \mu(R_i) |c_i|^2.
\]

We define the inner product on \( L^2(\Omega, \mu) \) in the standard form:

\[
\langle f, g \rangle = \int_\Omega f(x) g(x) \, d\mu(x),
\]

and thus the squared norm of \( f \) is:

\[
\|f\|^2 = \langle f, f \rangle = \int_\Omega |f(x)|^2 \, d\mu(x).
\]

Assuming the transform is complete (invertible) but not orthonormal, the coefficient energy \( E_f \) may not equal \( \|f\|^2 \), but it is bounded above and below by frame constants. Therefore, we establish the frame inequality\cite{daubechies1992ten}:

\[
A \|f\|^2 \leq \sum_i \mu(R_i) |c_i|^2 \leq B \|f\|^2,
\]

for frame bounds \( A, B > 0 \). These constants describe the stability and fidelity of the reconstruction.

In the case of fixed refinement parameters and a uniform transform structure, we define the total measure covered by the refined regions as:

\[
M := \sum_i \mu(R_i).
\]

We can then define a normalized frame inner product as:

\[
\langle f, f \rangle_{\text{frame}} := \frac{1}{M} \sum_i \mu(R_i) |c_i|^2,
\]

which approximates the true \( L^2 \) norm:

\[
\langle f, f \rangle_{\text{frame}} \approx \|f\|^2.
\]

This implies a natural relation:

\[
\sum_i \mu(R_i) |c_i|^2 = M \|f\|^2,
\]

so that the frame bounds may be expressed as:

\[
A = M, \quad B = 1 \quad \text{or} \quad A = 1, \quad B = M,
\]

depending on whether one takes the function norm or coefficient sum as the normalization standard.

It is worth reminding here that frame bounds obey the following:

\begin{theorem}[Frame Bounds]
There exist constants \( A, B > 0 \) such that for all \( f \in L^2(\Omega) \),
\[
A \|f\|_{L^2}^2 \leq \sum_{m,i} |c_{m,i}|^2 \cdot |\Omega_{m,i}| \leq B \|f\|_{L^2}^2.
\]
\end{theorem}

This\cite{daubechies1992ten} is a direct result of the basis not being orthonormal. Because each point of f is reconstructed using each region that contains that point in its indicator function - this implies lots of overlapping contributions. The lack of freedom to choose these coefficient values in multiple ways to reconstruct f is directly related to the fact that our transform is not overcomplete. 

This framework connects the measure-theoretic structure of the decomposition with the classical inner product, and establishes the precise role of the refinement geometry in determining transform stability.

Comment: The frame bounds follow a geometric series relating to the measure of regions at each order and index. For consistent self similar refinements this looks like a simple geometric series. For more complicated geometries this follows combinations of partial geometric sums.

\section{Connection to Existing Transforms}

\subsection{Radon Transform}

Lets consider first the Radon Transform in 2d. Set the refinement multiplicity to 2, dispersion to 1/2. Instead of taking an N-gon just take two points with the center point at the middle. If we set our angle to 0 (pointing right to left) and we set our dispersion to 1/2 the decomposition will break up $\Omega$ into up and down regions, as thin as we want. If we attempt to compute a reconstruction of the original function f this will introduce line artifacts.Let \( T_{\theta} \) be the reconstruction transform for \( f \) corresponding to an angle \( \theta \). We then have \cite{deans2007radon}

\begin{equation}
    f(x, y) = \lim_{N \to \infty} \frac{\sum_{i=1}^{N} T_{\frac{180 \cdot i}{N}} f(x, y)}{N} 
\end{equation}

Comments: With arbitrarily high fidelity data I believe we can relax this conditon on needing a limiting N case, how can we show this?

\subsection{Generalization to Higher Dimensions}

We now extend the concept of the Radon Transform to \( \mathbb{R}^N \). In this framework, the refinement multiplicity plays a critical role in determining the presence of reconstruction artifacts.

\subsection{The Radon Zone and Subspace Artifacts}

Let \( f \in BV(\Omega) \), where \( \Omega \subset \mathbb{R}^N \). Let \( T_{\Theta} f \) denote the Voronoi refinement transform applied along directional configuration \( \Theta \), with fixed refinement multiplicity \( M < N+1 \), dispersion \( \delta \in (0, 1) \), and hierarchical Voronoi-based partitioning.

We define the \textit{Radon Zone} as the regime in which the refinement multiplicity is strictly less than \( N+1 \). In this regime, reconstructions \( T_{\Theta} f \) inherently exhibit artifacts aligned with lower-dimensional subspaces of dimension \( N - M \).

\subsection{Generalized Geometric Refinement Radon Transform}

\begin{theorem}[Generalized Geometric Refinement Radon Transform]
Let \( f \in BV(\Omega) \), with \( \Omega \subset \mathbb{R}^N \). If the refinement multiplicity is \( M < N+1 \), then:
\begin{enumerate}
    \item Reconstruction using \( T_{\Theta} f \) without angular averaging exhibits artifacts of dimension \( N - M \).
    \item To suppress these artifacts and recover \( f \), one must average over configurations \( \{ \Theta_i \} \) drawn from rotational symmetries or optimal angular distributions.
    \item Specifically, if the artifact dimension is \( N - n \), then the averaging directions should be chosen as solutions to the generalized Thomson problem on \( S^{n} \subset S^{N-1} \).
\end{enumerate}
Define the averaged transform:
\[
f_K(x) = \frac{1}{K} \sum_{i=1}^{K} T_{\Theta_i} f(x)
\]
where \( \{ \Theta_i \} \) are direction sets distributed according to the optimal spherical configuration. Then:
\[
\lim_{K \to \infty} f_K(x) = f(x)
\]
in the appropriate norm (e.g., \( L^2(\Omega) \)) as both the angular sample count \( K \to \infty \) and refinement depth increase. \cite{deans2007radon}
\end{theorem}

This generalization captures a hierarchy of dimensional artifacts and highlights the geometric necessity of sufficient angular diversity in directional transforms for accurate recovery.

\section{Complete Refinement Zone} 

\subsection{Symmetry}

Let us now turn our attention away from redundant decompositions. In the vast parameter space it may seem intractible to sit through checking all of them. We ask the question: What Voronoi refinements lead to self-symmetric refinements. A keen reader will have noticed that we have shown a transform that  decomposes an arbitrary region but has no guarantees of self similarity across the parameter space. The key insight goes as follows:

\begin{proposition}[Symmetric Refinements] If the daughter points generate Voronoi cells such that each daughter point is at the centroid of the cell that it produces, then the refinement is symmetric. \cite{du1999centroidal} 
\end{proposition}

The "trick" then is to pick the right refinement characteristics that allow for highly symmetric Voronoi regions. We have left this static throughout the refinement for this paper, but this is by no means a requirement. It just helps the initial formulation be tractable.

\subsection{Wavelet Transform}

Let us consider a quad-adic decomposition. Let us set the angle to  \( \frac{\pi}{4} \), the dispersion to 1/2, and the refinement multiplicity to 4. This naturally decomposes a square grid into \( 4^m\) regions, where m is the refinement step \cite{mallat2008wavelet}.

\section{Completeness, Uniqueness, Invertibility, Closure, and Stability}
To reconstruct \( f \) at any point, we consider the set of all parent regions of that point and sum the corresponding coefficients at each level of refinement. 

Since the space of piecewise continuous functions \( \mathcal{P} \) is a subset of the locally integrable functions \( L^1_{\text{loc}}(\Omega) \), and both are dense in \( L^2(\Omega) \), we present the following theorems in the more general setting of \( L^1_{\text{loc}} \). This choice preserves the constructive intuition derived from \( \mathcal{P} \) while maximizing generality and applicability. As a result, the key properties of completeness, uniqueness, and invertibility extend naturally to all of \( L^2 \) via standard limiting arguments.

We prove here completeness, uniqueness, and invertibility. We first take a quick break for exposition clarity to remind of the following:

\subsection{Function Space Hierarchy}
We begin with the space of piecewise continuous functions \( \mathcal{P} \), which is dense in the space of functions of bounded variation \( \mathrm{BV}(\Omega) \), which is in turn dense in \( L^2(\Omega, \mu) \):
\[
\mathcal{P} \subset \mathrm{BV}(\Omega) \subset L^2(\Omega, \mu).
\]
Our transform is first defined and proven well-behaved on \( \mathcal{P} \), then extended to \( \mathrm{BV} \) via convergence, and finally to all of \( L^2 \) via density arguments.

\begin{theorem}[Completeness]
Let \( f \in L^1_{\text{loc}}(\Omega) \), and let \( \{ \Omega_{m, i} \} \) be the regions defined by the hierarchical Voronoi refinement with dispersion \( \delta \in (0, 1) \). Define the transform coefficients recursively as
\[
    c_{0,0} = \frac{1}{|\Omega|} \int_{\Omega} f(x)\, dx,
\]
\[
    c_{m, i} = \frac{1}{|\Omega_{m, i}|} \int_{\Omega_{m, i}} f(x)\, dx - \frac{1}{|\Omega_{m-1, j}|} \int_{\Omega_{m-1, j}} f(x)\, dx,
\]
where \( \Omega_{m,i} \subset \Omega_{m-1,j} \). Then the partial reconstructions
\[
    f_m(x) = \sum_{k=0}^m \sum_{i} c_{k,i} \cdot \chi_{\Omega_{k,i}}(x)
\]
where \( \chi_{\Omega_{k,i}}(x) \) is the indicator function for the region \( \Omega_{k,i} \) associated with that coefficient. This will converge to \( f(x) \) almost everywhere as \( m \to \infty \). That is,
\[
    \lim_{m \to \infty} f_m(x) = f(x) \quad \text{a.e. on } \Omega.
\]
This follows from the Lebesgue Differentiation Theorem, as the refinement regions shrink to points and the telescoping sum converges pointwise.
\end{theorem}

\begin{theorem}[Uniqueness]
If two functions $f$ and $g$ yield identical transform coefficients over every region in the refinement, then $f = g$ almost everywhere on $\Omega$.

\textit{Proof.} Let $h = f - g$. Then all transform coefficients for \( h \) vanish:
\[
c_{0,0}^h = \frac{1}{|\Omega|} \int_{\Omega} h(x) dx = 0,
\]
\[
c_{m, i}^h = \frac{1}{|\Omega_{m, i}|} \int_{\Omega_{m, i}} h(x)\, dx - \frac{1}{|\Omega_{m-1, j}|} \int_{\Omega_{m-1, j}} h(x)\, dx = 0.
\]
This implies the average of \( h \) over each cell equals the average over its parent. Since the base average is zero and every refinement preserves this, the average of \( h \) is zero on every region. As \( \text{diam}(\Omega_{m,i}) \to 0 \), the Lebesgue Differentiation Theorem implies \( h(x) = 0 \) almost everywhere.
\end{theorem}

\begin{theorem}[Invertibility]
Let $f \in L^1_{\text{loc}}(\Omega)$ be a function such that all transform coefficients vanish:
\[
c_{0,0} = \frac{1}{|\Omega|} \int_{\Omega} f(x)\, dx = 0,
\]
\[
c_{m, i} = \frac{1}{|\Omega_{m, i}|} \int_{\Omega_{m, i}} f(x)\, dx - \frac{1}{|\Omega_{m-1, j}|} \int_{\Omega_{m-1, j}} f(x)\, dx = 0.
\]
Then \( f = 0 \) almost everywhere on \( \Omega \).

\textit{Proof.} The vanishing of all coefficients implies that the average of \( f \) over each cell equals its parent, and the global average is zero. By induction over refinement levels, every region has zero average. As the regions shrink and form a covering of \( \Omega \), the Lebesgue Differentiation Theorem implies \( f(x) = 0 \) almost everywhere.
\end{theorem}

\begin{proposition}[Closure]
The transform space remains closed under refinement as $n \to \infty$, producing well-defined limiting hyperspherical sectors. Reducing the refinement multiplicity to 1 creates degeneracies but the transform remains well behaved.

The transform space remains well behaved but produces degenerate refinements at $\delta = 0$ or sometimes at $\delta = 1$. The points either stay at the center for $\delta = 0$ or go straight to the boundary at $\delta = 1$. This can cause further refinements to repeat earlier Voronoi diagrams, resulting in a well-behaved but incomplete transform. Rotations are continuous and well posed and do not admit any irregularity under limit investigation. All limits produce well behaved transforms.

Furthermore, closure is rigorously supported by the following inductive argument: for a given refinement multiplicity $N$, if a complete and non-degenerate refinement exists for all orders, and we can construct a valid refinement for $N+1$ while preserving non-degeneracy and coverage, then by induction, the transform space remains well defined for all refinement multiplicities $N \geq n_0$. This guarantees that arbitrarily high multiplicities yield consistent and stable transform behavior.
\end{proposition}

\section{Applications and Physical Interpretation}

\subsection{Entropy Dynamics and Microstate Interpretation}

One of the most intriguing implications of the transform lies in its ability to encode entropy-like behavior through the evolution of its coefficient spectrum. At each refinement level, the transform redistributes energy across finer regions, mimicking physical processes such as diffusion, dispersion, and energy cascade. This leads to a natural entropy formulation.

We define the coefficient spectrum \( \{ c_{m,i} \} \) at level \( m \), and associate to each region \( \Omega_{m,i} \) a normalized energy contribution:
\[
p_{m,i} = \frac{|c_{m,i}|^2 \cdot |\Omega_{m,i}|}{\sum_j |c_{m,j}|^2 \cdot |\Omega_{m,j}|}.
\]
Each \( p_{m,i} \) represents the proportion of the total signal energy contained within the region \( \Omega_{m,i} \). It captures both the magnitude of the coefficient and the measure of the region — reflecting the intuition that a high-energy contribution on a large domain is more significant than one on a small domain.

We then define a Shannon-like entropy at level \( m \) as:
\[
S_m = - \sum_i p_{m,i} \log p_{m,i},
\]
where \( \sum_i p_{m,i} = 1 \) by construction. As the refinement progresses, a single strong coefficient tends to split into many smaller contributions over space, increasing entropy.

This process closely parallels statistical mechanics: the transform's coefficients can be interpreted as emergent microstates of a coarse-to-fine energy distribution \cite{jaynes1957information}. In this analogy, the entropy functional behaves like the entropy of a canonical ensemble, and the measure-weighted energy contributions \( p_{m,i} \) resemble probabilities in the Gibbs distribution. In fact, one may define a formal partition function at level \( m \) as:
\[
Z_m = \sum_i \exp\left( -\beta E_{m,i} \right),
\]
where \( E_{m,i} = -\log p_{m,i} \) and \( \beta \) plays the role of an inverse temperature. This makes the entropy formula a direct analog of the thermodynamic identity:
\[
S = -\sum_i p_i \log p_i = \log Z + \beta \langle E \rangle.
\]

In physical terms, as energy disperses spatially (due to the refinement structure), the number of accessible microstates increases, and the system becomes more disordered — exactly in line with the Second Law of Thermodynamics. In the transform picture, this entropy growth corresponds to the decay of variance across the coefficient hierarchy.

We conjecture that this entropy functional, under suitable physical constraints on the refinement parameters, may serve as a variational principle for simulating energy-minimizing dynamics in physical systems. This opens the door to applications in thermodynamics, image compression, statistical inference, and even fundamental physics. 

In addition to this analysis there are multiple other existing fields that could benefit from this work:

\begin{itemize}
    \item \textbf{PDE Simulation:} Averaging over neighboring regions near centric refinements approximates PDE evolution via a discrete Pizzetti-like update. Centric refinements can use the center cell as a naturally surrounded region for the use of the Pizzetti formula. This allows for the updating of an input PDE according to the transform geometry selected.
    \item \textbf{Image and Signal Processing:} Adaptive refinement and hierarchical averaging allow scale-sensitive decomposition and reconstruction. This allows for correlation across scale in non-dyadic refinement schemes.
    \item \textbf{Nonlinear Systems:} The transform tracks structure formation in turbulence, gravity, and self-organizing systems. Every transform has a "natural" language that it uses to encode functional behavior. By leveraging different expected self similarity patterns in information preservation for black holes or for the genesis of an organism. 
    \item \textbf{Medical Imaging and Diagnostics:} Traditional imaging systems have relied on square grids due to mathematical constraints, but by leveraging a heirarchical transform theory it would seem reasonable to have an opportunity to deviate from this. This provides a new mechanism to observe fMRI images adapting in a maximally symmetric way to the brain, or perhaps even in a physically relevant way by allowing non-static refinements. 
    \item \textbf{General Image Processing} Besides the grid being generalizably adaptable, this has the advantage that wavelets currently have for enhancing edge characteristics. By plotting the high order coefficients the edge characteristics can be reconstucted very well which allows for a generalized automated segmentation strategy.
\end{itemize}

\section{Conclusion}
Our new geometric refinement-based transform provides a robust, mathematically complete generalization of classical transforms. By leveraging simplex-based refinement structures, the transform maintains completeness, uniqueness, invertibility, closure, and stability, naturally suited to self-similar functional structures. Given the plethora of natural phenomenon that exhibit self similarity it seems inevitable that this transform will help elucidate these patterns. Transform theory has provided (pun intended) transformative changes to a plethora of fields from signal processing to functional analysis to PDE analysis and medical imaging and therefore it seems appropriate that this transform would find an abundance in utility.

While the transform shows promise across a range of applications, further work is needed to optimize numerical implementation, explore convergence in non-Euclidean settings, and study operator theory in the transform domain. There is likely rich structure emerging from parameter-space correlations to manifolds and adaptive refinement schemes.

\bibliographystyle{plain}
\bibliography{references}

\end{document}